\documentclass{ifacconf}

\usepackage{graphicx}      
\usepackage{natbib}        
\usepackage{amsmath, amssymb}
\usepackage{mathrsfs}
\usepackage{epstopdf}

\newcommand{\bA}{{\bf A}}
\newcommand{\bB}{{\bf B}}

\newcommand{\bI}{{\bf I}}
\newcommand{\bK}{{\bf K}}
\newcommand{\bM}{{\bf M}}
\newcommand{\bx}{{\bf x}}

\newcommand{\bb}{{\bf b}}
\newcommand{\bc}{{\bf c}}

\newtheorem{theorem}{Theorem}[section]

\begin{document}
\begin{frontmatter}

\title{The AAA framework for modeling linear dynamical systems with quadratic output} 

\thanks[footnoteinfo]{The work of Gugercin was supported in parts by NSF through Grant DMS-1720257 and DMS-1819110.}
\thanks[footnoteinfo]{\textcopyright \ 2020 Ion Victor Gosea and Serkan Gugercin. This work has been accepted to IFAC for publication under a Creative Commons Licence CC-BY-NC-ND.}

\author[First]{Ion Victor Gosea} 
\author[Second]{Serkan Gugercin}

\address[First]{Data-Driven System Reduction and Identification Group, Max Planck Institute for Dynamics of Complex Technical Systems, Magdeburg, Germany (e-mail: gosea@mpi-magdeburg.mpg.de)}
\address[Second]{Department of Mathematics and Division of Computational Modeling and Data Analytics, Virginia Tech, Blacksburg, VA 24061, Virginia, USA (e-mail: gugercin@vt.edu)}

\begin{abstract}  
	We consider linear dynamical systems with quadratic output.
	We first define the two transfer functions, a single-variable and a multivariate one, that fully describe the dynamics of these special nonlinear systems. Then, using the samples of these two transfer functions, we extend the \textsf{AAA} algorithm to model linear systems with quadratic output from data.
\end{abstract}

\begin{keyword}
Model reduction, complex systems, time-invariant systems.
\end{keyword}

\end{frontmatter}

\section{Introduction}

Model order reduction (MOR) typically refers to a class of methodologies that can be used to approximate large-scale dynamical systems with much smaller systems that ideally have similar response characteristics as the original ones. Many MOR methods have been developed in the last decades; we refer the reader to  \cite{ACA05,BBF14,AntBG20,quarteroni2015reduced,siammorbook2017} and to the references therein for more details on various different approaches to model reduction. The approach we consider here falls under the systems theoretical (input/output) framework for model reduction.  

Many time-dependent processes modeled and studied in real-world applications exhibit nonlinear dynamics. In some cases, in order to simplify the modeling part, linearization of the dynamics is performed in many engineering branches. This is usually done around an operating point in a given domain and is restricted to local conditions. In order to obtain more general models, nonlinearities need not be omitted and, instead, need to be suitably treated. 
Hence, the study, analysis, and modeling of nonlinear dynamical systems have been the focus of considerable research in the last decades. Consequently, many model reduction methods that can be directly applied to nonlinear systems have been proposed. For a fairly extensive and recent review on MOR for nonlinear and linear systems, we refer the reader to \cite{BBF14}.

Considerable progress has been made by extending systems theoretic MOR methods from linear to certain classes of  nonlinear systems; such as bilinear systems, quadratic-bilinear (\textsf{QB}) systems, and 
linear systems with quadratic output (\textsf{LQO}), which is the main focus of our work; see, e.g.,  \cite{damm11,AGIbilinear,breiten12,flagg15} for extension of systems theoretical techniques to bilinear systems; 
\cite{gosea18,benner2018h2,morBTQBgoyal,kramer2019balanced} to
\textsf{QB} systems, and
\cite{PN19,BGP19}
 to \textsf{LQO} systems. Additional contributions to MOR of nonlinear systems were made by \cite{kawano2007} and by \cite{astolfi2010}.
 
 Moreover, system identification of nonlinear systems has been a popular topic for decades already. In particular, we mention here the special case of identifying linear systems with nonlinear output or input functions, e.g. the so-called Wiener and, respectively, Hammerstein models. Note also that LQO systems are a special class of Wiener models for which the nonlinear output mapping is quadratic. Significantly effort has been allocated for identification of such models; see, for example, \cite{juditsky95}, \cite{giri10} and the references therein.

We are interested data-driven methods for model reduction where one does not necessarily have access to internal degrees of freedom, i.e., the underlying large-scale state-space representation is unknown. Instead, one has access only to a collection of input/output measurements corresponding to the original system. In this work, by data we mean
frequency domain samples  
 of the input/output mapping, known as the transfer function, of the underlying system. Hence, we focus here on data-driven methods that use rational interpolation and least-square approaches to fit a rational function to the given set of frequency-domain (transfer function) measurements; such as the vector fitting approach, \cite{GS99}, the Loewner framework in \cite{mayo2007fsg,ALI17}, and the \textsf{AAA} algorithm in \cite{NST18}. This contribution aims at extending the \textsf{AAA} algorithm (Section \ref{sec:aaa}) to a special class of nonlinear systems outlined
 (Sections \ref{sec:lqo} and \ref{sec:prop}).

\section{The AAA algorithm} \label{sec:aaa}

Given a set of measurements $h_k = H(\xi_k)$ of an underlying (transfer) function $H(\cdot)$ for $k=1,2,\ldots,N_s$  where
$\xi_k \in \mathbb{C}$ are the \emph{sampling (support) points} and 
$h_k\in \mathbb{C}$ are the sampled data,
the \textsf{AAA} (Adaptive Antoulas-Anderson) algorithm computes a
rational function $r(s)$, either to a specified accuracy or of a specified order, that approximates the sampled data. The approximating rational function $r(s)$ is written in barycentric format, a numerically stable representation of rational functions:
\begin{equation}\label{eq:bary_rep}
r(s) =  \frac{\displaystyle\sum_{k=0}^n \frac{w_k h_k}{s - \xi_k}}{\displaystyle\sum_{k=0}^n\frac{w_k}{s - \xi_k}},
\end{equation}
where the \emph{weights} $w_k \in \mathbb{C}$ are to be determined.  It directly follows from the form of the rational function $r(s)$  in  \eqref{eq:bary_rep}  satisfies the  interpolation conditions $r(\xi_k) = h_k$ for $1 \leqslant k \leqslant n$, assuming $w_k \neq 0$. Hence, interpolation at the $n$ support points is attained for free by means of the special representation chosen. Then, there is a freedom in choosing the weights $w_k \neq 0$ to match the remaining the $N_s-n$ data points in an appropriate measure. 

Assuming enough degrees of freedom, the framework of  
\cite{AA86} chooses the weights $w_k$ to enforce interpolation at the remaining points as well. In contrast, the \textsf{AAA} algorithm is iterative, and combines interpolation and least-squares fitting. Assume at  step $k$,  a rational function $r_k(s)$ in barycentric form \eqref{eq:bary_rep}
is constructed to interpolate $k$ data points.
The idea is that the next support (interpolation) point $\xi_{k+1}$ is selected by means of a greedy algorithm, i.e., the current approximation error of $r_{k}(s)$ at $\xi_{k+1}$ is maximum within the data set. Afterwards, the corresponding weights $w_k$ are computed by solving a least-squares problem to minimize the least-square deviation in the remaining data points.
 For further details, we refer the reader to the original source \cite{NST18}.  The \textsf{AAA} algorithm proved very flexible and effective, and has been employed in various applications such as  rational approximation over disconnected domains \cite{NST18}, 
solving nonlinear eigenvalue problems \cite{lietaert2018automatic}, modeling of parametrized dynamics \cite{CS20}, and approximation of matrix-valued functions \cite{GG20}.

\section{Linear systems with quadratic output}
\label{sec:lqo} 

The main focus of this work is to study data-driven MOR of linear systems with quadratic output, which are described in state-space representation as
\begin{align}\label{eq:def_LQO}
\Sigma_{\rm LQO}: \begin{cases}  \dot\bx(t)=\bA\bx(t)+\bb u(t),\\
\hspace{0.5mm} y(t)=\bc^T\bx(t)+\underbrace{\bK \big{[} \bx(t) \otimes \bx(t) \big{]}}_{\bx^T(t) \bM \bx(t)}, \end{cases}
\end{align}
where  $\bA \in \mathbb{R}^{n \times n}$, $\bb, \bc \in \mathbb{R}^{n}$, $\bK \in \mathbb{R}^{1 \times n^2}$, and $ \bM \in \mathbb{R}^{n \times n} \ (\bK = \text{vec}(\bM))$.
Additionally, note that in (\ref{eq:def_LQO}), the symbol $\otimes$ denotes the Kronecker product of the vector $\bx = [\bx_1 \ \bx_2 \ \cdots \ \bx_n]^T$ with itself, i.e.
$$
\bx \otimes \bx = [\bx_1^2 \ \ \bx_1 \bx_2 \ \ \bx_1 \bx_3 \ \ \cdots \ \ \bx_1 \bx_{n} \ \ \cdots \bx_{n}^2]^T \in \mathbb{R}^{n^2}.
$$
Several MOR methodologies have been already proposed for dealing with the case of \textsf{LQO} systems. More precisely, balanced truncation-type methods were considered in \cite{BM10,PN19,BGP19}, while interpolation-based methods were used in \cite{BNLM12,GA19}.
In some cases, the vector $\bc^T$ in (\ref{eq:def_LQO}) is $\bc = \mathbf{0}$, and thus the output has only the quadratic term.

For \textsf{LQO} systems in (\ref{eq:def_LQO}), the nonlinearity is present in the state-to-output equation only. Hence, one can write the input-output mapping of system $\Sigma_{\rm LQO}$ in the frequency domain using two transfer functions:
\begin{enumerate}
\item one corresponding to the linear part of the output, i.e., $y_1(t) = \bc^T \bx(t)$;
\item one corresponding to the quadratic part  of the output, i.e., $y_2(t)  = \bK (\bx(t) \otimes \bx(t))$.
\end{enumerate}

The linear transfer function corresponding to $y_1(t)$, denoted by $H_1(s)$, is as in the classical case of linear systems. It is a rational function and can be written in terms of the system matrices as 
\begin{equation} \label{eq:H1}
    H_1(s) = \bc^T(s\bI_n-\bA)^{-1} \bb,
\end{equation}
where $\bI_n$ is the identity matrix of dimension $n \times n$. Additionally, the quadratic transfer function, corresponding to the output $y_2(t)$ is also a rational function, but of two variables. It is defined as follows 
\begin{equation} \label{eq:H2}
H_2(s,z) = \bK \Big[ (s\bI_n-\bA)^{-1} \bB \otimes (z\bI_n-\bA)^{-1} \bb \Big].
\end{equation} 
We will use an extension of the barycentric representation (\ref{eq:bary_rep}) to represent  $H_2(s,z)$ in a suitable barycentric-like format to extend \textsf{AAA} to \textsf{LQO} systems in \eqref{eq:def_LQO}. 

\section{Proposed framework for data-driven modeling of \textsf{LQO} systems}
\label{sec:prop}
Recall that the \textsf{AAA} framework for linear systems uses (transfer) function samples $h_k = H(\xi_k)$. 
In extending \textsf{AAA} to \textsf{LQO} systems, data (measurements) will correspond to sampling not only the single-variable transfer function $H_1(s)$ in \eqref{eq:H1} but also the two-variable transfer function 
$H_2(s,z)$ in \eqref{eq:H2}. 

\begin{theorem}\label{theo1}
The barycentric form for approximating the first transfer function
$H_1(s)$ is denoted with $r_1(s)$ and can be written in a similar manner as \eqref{eq:bary_rep}, namely
\begin{equation}\label{eq:H1_bary}
r_1(s) = \frac{\displaystyle \sum_{k=1}^n \frac{ w_k h_k }{s - \xi_k }}{1+\displaystyle \sum_{k=1}^n \frac{w_k }{s- \xi_k}}.
\end{equation}
It follows that $r_1(s)$ interpolates the data at the support points $\{\xi_1, \xi_2, \ldots, \xi_n\}$, i.e.,
\begin{equation}
r_1(\xi_k) = w_i = H_1(\xi_k), \ (\forall) \ i=1,2,\ldots,n.
\end{equation}
\end{theorem}
Note that the scalar $1$ was added to the denominator to make the rational approximant $r_1(s)$ strictly proper. However, this does not affect the interpolation property.

In what follows we will sketch the proof of Theorem \ref{theo1}. Define the polynomial $P(s) = \prod_{k=1}^n (s-\xi_k)$ and also $P_i(s) = \prod_{k=1, k \neq i}^n (s-\xi_k) \ \forall \ i=1,\ldots,n$. Then, by bringing all the fraction terms of $r_1(s)$ in (\ref{eq:H1_bary}) to the same common denominator $P$, one can rewrite the rational function as
\begin{align}
\begin{split}
r_1(s) = \frac{\sum_{k=1}^n w_k h_k P_k(s)}{P(s) + \sum_{k=1}^n w_k P_k(s)} 
\Rightarrow r_1(\xi_i) = \frac{w_i h_i P_i(\xi_i)}{w_i P_i(\xi_i)} = h_i
\end{split}
\end{align}
Next, we need to pick an appropriate barycentric-like form for approximating $H_2(s,z)$. The form for $H_2(s,z)$ cannot be chosen independently from that of $H_1(s)$ (e.g., as a two-variable barycentric form as in \cite{IA14}) because  $H_1(s)$ and $H_2(s,z)$ share the same poles and the weights are related. This is easy to see from the transfer function $H_1(s)$ and $H_2(s,z)$. The eigenvalues of the same matrix $\bA$ contribute to the poles of both $H_1(s)$ and $H_2(s,z)$. In other words, we need to pick the barycentric(-like) forms for $r_1(s)$ and $r_2(s,z)$ to guarantee that they correspond to an underlying \textsf{LQO} system.
We have the following result:

\begin{theorem} \label{thm:lqobar}
Given the interpolation points
$\xi_k$ for $k=1,2,\ldots,n$, let 
the data 
$h_k = H_1(\xi_k)$ and
$h_{k,\ell} = H_2(\xi_k,\xi_l)$
for $k,\ell = 1,\ldots,n$
result form sampling the transfer functions $H_1(s)$ and $H_2(s,z)$ corresponding to 
\textsf{LQO} system in
\eqref{eq:def_LQO}. Then, 
$r_1(s)$ in \eqref{eq:H1_bary} interpolates the data $h_k$, i.e., $r_1(\xi_k) = H_1(\xi_k)$ for $k=1,2,\ldots,n$.
Define, the two-variable  function $r_2(s,z)$ in a barycentric-like form:
\scriptsize
\begin{equation}\label{eq:H2_bary}
    r_2(s,z) = \frac{\displaystyle \sum_{k=1}^n \sum_{\ell=1}^n \frac{ h_{k,\ell}w_k w_\ell}{(s-\xi_k)(z-\xi_\ell)}}{1+\displaystyle \sum_{k=1}^n \frac{w_k }{s- \xi_k} + \displaystyle \sum_{\ell=1}^n \frac{w_\ell }{z- \xi_\ell} +\displaystyle \sum_{k=1}^n \sum_{\ell=1}^n \frac{ w_k w_\ell}{(s-\xi_k)(z-\xi_\ell)}},
\end{equation}
\normalsize
where $h_{k,\ell} := H_2(\xi_k,\xi_\ell)$. Then, $r_2(s,z)$ interpolates the data  $h_{k,\ell}$, i.e.,
$r_2(\xi_k,\xi_\ell) =
h_{k,\ell} = H_2(\xi_k,\xi_\ell) $ 
for $k,\ell = 1,\ldots,n$. Moreover, $r_1(s)$ and $r_2(s,z)$ correspond to a reduced \textsf{LQO} model. In others words, there exist a \textsf{LQO} model as in \eqref{eq:def_LQO} whose first (linear) transfer function is $r_1(s)$ and second transfer function is $r_2(s,z)$. 
\end{theorem}
The interpolation property stated in Theorem \ref{thm:lqobar} can be proven in a similar manner to Theorem \ref{theo1}, using polynomials in two variables. We skip those details here.

Note that the barycentric representations for the rational interpolants $r_1(s)$ and $r_2(s,z)$ have common weights given by variables $w_k$, and by products $w_k w_\ell$, respectively. This indeed shows that the two rational functions are connected to each other. Moreover, Theorem \ref{thm:lqobar} reveals how to construct a reduced \textsf{LQO} system directly from  transfer function samples. And the corresponding barycentric forms for $r_1(s)$ and $r_2(s,z)$
directly lend themselves to extending
\textsf{AAA} to modeling \textsf{LQO} systems. The forms 
$r_1(s)$ and $r_2(s,z)$ interpolate the data by construction.

 As in the original \textsf{AAA} algorithm for linear dynamics, we will produce the approximants iteratively. More precisely, we increase the orders of $r_1(s)$
and $r_2(s,z)$ and will automatically interpolate a subset of the data. Then, the method will choose the free parameters given by the weights $w_k$ in order to minimize the least-squares distance in the remaining data.

\bibliography{LQO_AAA_ref}             

\end{document}